\documentclass{amsart}
\usepackage[dvips,line,matrix,frame,curve,poly]{xy}
\usepackage{amsmath,xspace,url, textcomp}
\newtheorem{thm}{Theorem}[section]
\newtheorem{prop}[thm]{Proposition}
\newtheorem{dfn}[thm]{Definition}

\newcommand{\Dfn}[1]{\emph{#1}}                            

\author{Martin Rubey}
\thanks{Research partially supported by the Austrian
    Science Foundation FWF, grant S9607-N13, in the framework of the National
    Research Network \lq\lq Analytic Combinatorics and Probabilistic Number
    Theory\rq\rq.}

\address{\url{martin.rubey@univie.ac.at}, Universit\"at Wien}

\keywords{nestings and crossings of matchings and permutations, PDSAWs}
\title{Nestings of Matchings and Permutations and North Steps in PDSAWs}
\usepackage[breaklinks=true,hyperref]{hyperref}
\begin{document}
\begin{abstract}
  We present a simple bijective proof of the fact that matchings of $[2n]$ with
  $N$ nestings are equinumerous to \Dfn{partially directed self avoiding walks}
  confined to the symmetric wedge defined by $y=\pm x$, with $n$ east steps and
  $N$ north steps.  A very similar construction connects permutations with $N$
  nestings and \Dfn{PDSAWs} remaining below the $x$-axis, again with $N$ north
  steps.  Furthermore, both bijections transport several combinatorially
  meaningful parameters.
\end{abstract}
\maketitle
\section{Introduction}
\label{sec:introduction}

This article exhibits a connection, at first maybe surprising, between two at
present very actively researched, yet classic, areas of combinatorics.  The
first area to be mentioned concerns the enumeration of matchings,
set-partitions and permutations, keeping track of various statistics such as
crossings and nestings.  For example, it was observed only rather recently by
Martin Klazar and Marc Noy (see Section~1 of~\cite{Klazar2006}) that the joint
distribution of crossings and nestings of matchings is symmetric.  A little
later, this proved to be true also for set-partitions, as shown by Anisse
Kasraoui and Jiang Zeng~\cite{KasraouiZeng2006}, and finally by Sylvie
Corteel~\cite{Corteel2007} for permutations.

The other area concerns counting self avoiding walks on lattices under various
restrictions.  These objects are not only interesting from a purely
combinatorial point of view, but also for physicists, who seem to use them as
models for polymers in dilute solutions.  Unfortunately, self avoiding walk
models are usually intractable from a combinatorial point of view.  However,
imposing some sort of directedness on the walks, we obtain models that are
easier to deal with.  More precisely, in this article we will consider self
avoiding walks that are \Dfn{partially directed}:
\begin{dfn}
  A \Dfn{partially directed self avoiding walk}, short \Dfn{PDSAW} is a walk in
  the plane, starting at the origin and taking unit east, north, and south
  steps, where, however, a north step must not be immediately followed or
  preceded by a south step.  Within this article, we restrict our attention to
  walks that either stay within the symmetric wedge, defined by $y=\pm x$, or
  the asymmetric wedge enclosed by the $x$-axis and $y=-x$.  Furthermore, we
  require that the paths end at some point on $y=-x$.
\end{dfn}

\subsection{PDSAWs in the symmetric wedge}
\label{sec:symmetric-intro}

\begin{figure}
  \begin{tabular}{cc}
    \begin{xy}<11pt,0pt>:
      {\def\n{12}},
      (13,0);(25,\n);**[|(3)]@{-},
      (13,0);(25,-\n);**[|(3)]@{-},
      @i @={(1,0),(0,-1),(1,0),(0,3),(1,0),(0,0),(1,0),(0,-2),(1,0),(0,2),(1,0),(0,-2),(1,0),(0,-7),(1,0),(0,1),
        (1,0),(0,-3),(1,0),(0,-1),(1,0),(0,5),(1,0),(0,-7)},
      (13,0)="prev", @@{;p+"prev";"prev";**[|(3)]@{-}="prev"},
      (20,-7);(25,-2);**[|(1.3)]@{.},
    \end{xy}
    &
    \begin{xy}<11pt,0pt>:
      {\def\n{14}},
      (15,0.1);(29,0.1);**[|(3)]@{-},
      (15,0);(29,-\n);**[|(3)]@{-},
      @i @={(1,0),(0,0),(1,0),(0,-2),(1,0),(0,0),(1,0),(0,0),(1,0),(0,1),(1,0),(0,-5),(1,0),(0,-1),(1,0),(0,-1),
        (1,0),(0,1),(1,0),(0,0),(1,0),(0,-2),(1,0),(0,2),(1,0),(0,0),(1,0),(0,-7)},
      (15,0)="prev", @@{;p+"prev";"prev";**[|(3)]@{-}="prev"},
      (21,-5.9);(29,-5.9);**[|(1.3)]@{.},
      (22,-6.9);(29,-6.9);**[|(1.3)]@{.},
    \end{xy}\\
    a. A PDSAW in the symmetric wedge & b. A PDSAW in the asymmetric wedge \\
    with $12$ east and $11$ north steps, & with $14$ east and $4$ north steps,\\
    $2$ factors and length of last descent $7$. & $3$ factors and length of
    last descent $7$.
  \end{tabular}
  \caption{Examples for PDSAWs in wedges.}
  \label{fig:PDSAW}
\end{figure}
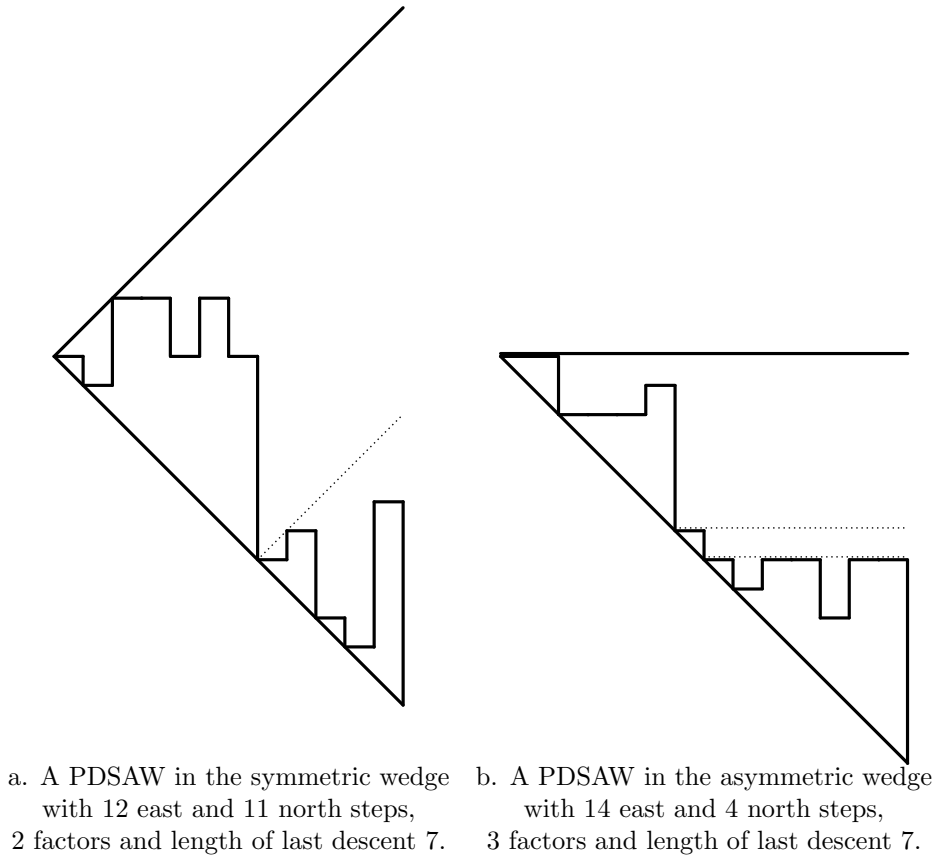

An example of a walk in the symmetric wedge can be found in
Figure~\ref{fig:PDSAW}.a.  Since a walk is entirely determined by the
$y$-coordinates of its east steps, it is immediate that the total number of
walks with $n$ east steps is given by $(2n-1)!!=1\cdot 3\dots (2n-3)(2n-1)$.
Note that this is also the number of matchings of the set
$[2n]=\{1,2,\dots,2n\}$.

By introducing the iterated kernel method, Janse van Rensburg, Thomas Prellberg
and Andrew Rechnitzer~\cite{RensburgPrellbergRechnitzer2007} were able to
derive a rather complicated expression for the generating function of PDSAWs
according to the \emph{total} number of steps.  Using that expression they also
computed the asymptotic number of such paths.

Intriguingly, the generating function they found was composed of pieces which
all seemed to have a combinatorial interpretation.  Roughly, it consisted of an
alternating series of powers of Catalan generating functions, albeit shifted by
a quadratically growing power, hinting at the distinct possibility of a more
direct combinatorial derivation.  When one of them related this to Philippe
Flajolet, he pointed out an apparent similarity to a formula counting the
number of matchings with respect to crossings.  Following up on this
connection, they discovered the surprising fact that the generating function
for PDSAWs with $n$ east steps, where $q$ marks the number of north steps, is
as simple as

\begin{equation}
  \label{eq:Touchard}
  M_{2n}=\frac{1}{(1-q)^n}
  \sum_{i\geq0}(-1)^i\left(\binom{2n}{n-i}-\binom{2n}{n-i-1}\right)
  q^{\binom{i+1}{2}}.
\end{equation}

This formula is known as \Dfn{Touchard-Riordan Formula}, and counts the number
of matchings of the set $[2n]$ according to crossings, i.e., pairs of matched
points $\{i,j\}$ and $\{k,l\}$ with $i<k<j<l$, pictorially:
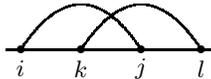
\begin{center}
{\setlength{\unitlength}{1mm}
\begin{picture}(28,10)(-30,-3)
\put(-30,0){\line(1,0){28}}
\put(-28,0){\circle*{1}}\put(-28,0){\makebox(0,-5)[c]{\small $i$}}
\put(-20,0){\circle*{1}}\put(-20,0){\makebox(0,-5)[c]{\small $k$}}
\put(-12,0){\circle*{1}}\put(-12,0){\makebox(0,-5)[c]{\small $j$}}
\put(-4,0){\circle*{1}}\put(-4,0){\makebox(0,-5)[c]{\small $l$}}
\qbezier(-28,0)(-20,12)(-12,0) \qbezier(-20,0)(-12,12)(-4,0)
\end{picture}}
\end{center}
A bijective proof of the latter fact was given by Jean-Guy
Penaud~\cite{Penaud1995}, exhibiting a whole zoo of combinatorial objects
counted by these numbers, including certain horizontally convex polyominoes.

Of course, we would now like to see a bijection between PDSAWs in the symmetric
wedge with $N$ north steps and matchings with $N$ crossings.  It is the main
purpose of this article to describe such a bijection.  However, it turns out
that it is more convenient to consider matchings with $N$ nestings instead,
that is, pairs of matched points $\{i,j\}$ and $\{k,l\}$ with $i<k<l<j$:
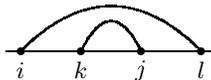
\begin{center}
{\setlength{\unitlength}{1mm}
\begin{picture}(28,10)(-30,-3)
\put(-30,0){\line(1,0){28}}
\put(-28,0){\circle*{1}}\put(-28,0){\makebox(0,-5)[c]{\small $i$}}
\put(-20,0){\circle*{1}}\put(-20,0){\makebox(0,-5)[c]{\small $k$}}
\put(-12,0){\circle*{1}}\put(-12,0){\makebox(0,-5)[c]{\small $j$}}
\put(-4,0){\circle*{1}}\put(-4,0){\makebox(0,-5)[c]{\small $l$}}
\qbezier(-28,0)(-16,12)(-4,0) \qbezier(-20,0)(-16,8)(-12,0)
\end{picture}}
\end{center}
Note that nestings and crossings in matchings (and, in fact, also in partitions
and permutations, given appropriate definitions) are equidistributed, see the
articles by Anisse Kasraoui and Jiang Zeng, and Sylvie
Corteel~\cite{Corteel2007, KasraouiZeng2006}.

The bijection we present in Sections~\ref{sec:KasraouiZeng}
and~\ref{sec:bijection} will transport several combinatorially meaningful
statistics on PDSAWs and matchings.  In particular, it will map factors of
PDSAWs to factors of matchings:
\begin{dfn}
  A \Dfn{factor} of a PDSAW in the symmetric wedge is a sub-path starting at
  $(a,-a)$ and ending at $(b,-b)$, such that
  \begin{itemize}
  \item[] all east steps after the point $(a,-a)$ are below the
    line $y=x-2a$,
  \item[] all east steps after the point $(b,-b)$ are below the line $y=x-2b$
  \end{itemize}
  for some $a$ and some $b$.

  A \Dfn{factor} of a matching is a sub-matching such that all elements of the
  sub-interval $\{a,a+1,\dots,b\}$ are matched within that interval.
\end{dfn}
As an example, the PDSAW in Figure~\ref{fig:PDSAW}.a has two (prime) factors
which are separated by the dotted line.

\begin{figure}
\begin{center}
{\setlength{\unitlength}{4mm}
\begin{picture}(24,3.6)(0.5,-0.7)
\qbezier(1,0)(4.5,4.0)(8,0)\qbezier(2,0)(6.0,4.5)(10,0)\qbezier(3,0)(3.5,1.0)(4,0)
\qbezier(5,0)(7.0,2.5)(9,0)\qbezier(6,0)(6.5,1.0)(7,0)\qbezier(11,0)(14.5,4.0)(18,0)
\qbezier(12,0)(16.5,5.0)(21,0)\qbezier(13,0)(13.5,1.0)(14,0)\qbezier(15,0)(17.0,2.5)(19,0)
\qbezier(16,0)(19.5,4.0)(23,0)\qbezier(17,0)(18.5,2.0)(20,0)\qbezier(22,0)(23.0,1.5)(24,0)
\put(0.5,0){\line(1,0){24}}
\put(1,0) {\circle*{0.2}}\put(1,0) {\makebox(0,-1)[c]{\small $1$}} 
\put(2,0) {\circle*{0.2}}\put(2,0) {\makebox(0,-1)[c]{\small $2$}} 
\put(3,0) {\circle*{0.2}}\put(3,0) {\makebox(0,-1)[c]{\small $3$}} 
\put(4,0) {\circle*{0.2}}\put(4,0) {\makebox(0,-1)[c]{\small $4$}} 
\put(5,0) {\circle*{0.2}}\put(5,0) {\makebox(0,-1)[c]{\small $5$}} 
\put(6,0) {\circle*{0.2}}\put(6,0) {\makebox(0,-1)[c]{\small $6$}} 
\put(7,0) {\circle*{0.2}}\put(7,0) {\makebox(0,-1)[c]{\small $7$}} 
\put(8,0) {\circle*{0.2}}\put(8,0) {\makebox(0,-1)[c]{\small $8$}} 
\put(9,0) {\circle*{0.2}}\put(9,0) {\makebox(0,-1)[c]{\small $9$}} 
\put(10,0){\circle*{0.2}}\put(10,0){\makebox(0,-1)[c]{\small $10$}}
\put(11,0){\circle*{0.2}}\put(11,0){\makebox(0,-1)[c]{\small $11$}}
\put(12,0){\circle*{0.2}}\put(12,0){\makebox(0,-1)[c]{\small $12$}}
\put(13,0){\circle*{0.2}}\put(13,0){\makebox(0,-1)[c]{\small $13$}}
\put(14,0){\circle*{0.2}}\put(14,0){\makebox(0,-1)[c]{\small $14$}}
\put(15,0){\circle*{0.2}}\put(15,0){\makebox(0,-1)[c]{\small $15$}}
\put(16,0){\circle*{0.2}}\put(16,0){\makebox(0,-1)[c]{\small $16$}}
\put(17,0){\circle*{0.2}}\put(17,0){\makebox(0,-1)[c]{\small $17$}}
\put(18,0){\circle*{0.2}}\put(18,0){\makebox(0,-1)[c]{\small $18$}}
\put(19,0){\circle*{0.2}}\put(19,0){\makebox(0,-1)[c]{\small $19$}}
\put(20,0){\circle*{0.2}}\put(20,0){\makebox(0,-1)[c]{\small $20$}}
\put(21,0){\circle*{0.2}}\put(21,0){\makebox(0,-1)[c]{\small $21$}}
\put(22,0){\circle*{0.2}}\put(22,0){\makebox(0,-1)[c]{\small $22$}}
\put(23,0){\circle*{0.2}}\put(23,0){\makebox(0,-1)[c]{\small $23$}}
\put(24,0){\circle*{0.2}}\put(24,0){\makebox(0,-1)[c]{\small $24$}}
\end{picture}}
\end{center}
\caption{The matching of $[24]$ corresponding to the PDSAW in
  Figure~\ref{fig:PDSAW}.a.}
  \label{fig:Matching}
\end{figure}
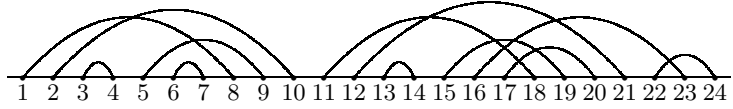

\begin{thm}\label{thm:pdsaws}
  PDSAWs with $n$ east steps,
  \begin{itemize}
  \item $N$ north steps,
  \item parity of area enclosed by the path and the line $y=-x$ equal to $A$,
    and
  \item length of the last descent equal to $M-1$
  \end{itemize}
  are in bijection with matchings of the set $[2n]$ with
  \begin{itemize}
  \item $N$ nestings,
  \item parity of the number of crossings equal to $A$, and
  \item $1$ being matched with $M$.
  \end{itemize}
  Moreover, factors are preserved, the last factor of the PDSAW being mapped to
  the first factor of the matching.
\end{thm}

(To compute the area enclosed by the path and the line $y=-x$, we only count
full squares and disregard triangles. The length of the last descent is the
number of south steps after the last east step of the PDSAW.)

\subsection{PDSAWs in the asymmetric wedge}
\label{sec:asymmetric-intro}

We now shift our attention to PDSAWs in the asymmetric wedge.  Again, it is
immediate that the total number of such walks with $n$ east steps equals $n!$.
Indeed, after performing a few computations and consulting the on-line
encyclopedia of integer sequences~\cite{OEIS}, one will be convinced that these
paths should be in bijection with permutations, north steps being mapped to
nestings as defined by Sylvie Corteel~\cite{Corteel2007}, or, alternatively, to
the number of occurrences of the (generalised) pattern $31\mbox{-}2$.  Similar
to the case of matchings, there is a formula for the generating function of
permutations of $[n]$ according to nestings, as shown by Sylvie
Corteel~\cite{Corteel2007} building on work of Lauren Williams~\cite{MR2102660}:

\begin{equation}
  \label{eq:Williams}
  P_n=\sum_{k=1}^n q^{-k^2}\sum_{i=0}^{k-1}(-1)^i [k-i]_q^n q^{ki}
  \left(\binom{n}{i}q^{k-i}+\binom{n}{i-1}\right)
\end{equation}

Sure enough, just a few days (or were it mere hours?)\ after being presented
with this conjecture, Philippe Nadeau came up with a surprisingly simple
bijection between PDSAWs that remain below the $x$-axis with $N$ north steps
and permutations that contain the generalised pattern $31\mbox{-}2$ $N$ times.
Motivated by this, and given the definitions of Sylvie
Corteel~\cite{Corteel2007}, it was not hard to find another bijection to
permutations with $N$ nestings, that again transports factors and the length of
the last descent nicely.  Of course, in this situation, a \Dfn{factor} is
defined slightly differently:

\begin{dfn}
  A \Dfn{factor} of a PDSAW in the asymmetric wedge is a sub-path starting at
  $(a,-a)$ and ending at $(b,-b)$, that stays below the line $y=-a$, for some
  $a$ and some $b$.

  A \Dfn{factor} of a permutation is a sub-permutation such that all elements
  less than $a$ are mapped to elements less than $a$, and all elements greater
  than $b$ are mapped to elements greater than $b$.
\end{dfn}

\newsavebox{\perm}
\sbox{\perm}{$\left(\begin{smallmatrix}
      1&2&3&4&5&6&7&8&9&10&11&12&13&14\\
      7&2&5&1&4&3&6&8&13&9&10&14&12&11
    \end{smallmatrix}\right)$}
\begin{figure}
\begin{center}
{\setlength{\unitlength}{4mm}
\begin{picture}(14,4)(0.5,-2)
\put(2,0.3){\oval(0.3,0.6)}
\put(8,0.3){\oval(0.3,0.6)}
\qbezier(1,0)(4.0,3.5)(7,0)\qbezier(3,0)(4.0,1.5)(5,0)\qbezier(9,0)(11.0,2.5)(13,0)
\qbezier(12,0)(13.0,1.5)(14,0)
\qbezier(4,-1)(2.5,-3.0)(1,-1)\qbezier(5,-1)(4.5,-2.0)(4,-1)\qbezier(6,-1)(4.5,-3.0)(3,-1)
\qbezier(7,-1)(6.5,-2.0)(6,-1)\qbezier(10,-1)(9.5,-2.0)(9,-1)\qbezier(11,-1)(10.5,-2.0)(10,-1)
\qbezier(13,-1)(12.5,-2.0)(12,-1)\qbezier(14,-1)(12.5,-3.0)(11,-1)
\put(1,0) {\makebox(0,-1)[c]{\small $1$}} 
\put(2,0) {\makebox(0,-1)[c]{\small $2$}} 
\put(3,0) {\makebox(0,-1)[c]{\small $3$}} 
\put(4,0) {\makebox(0,-1)[c]{\small $4$}} 
\put(5,0) {\makebox(0,-1)[c]{\small $5$}} 
\put(6,0) {\makebox(0,-1)[c]{\small $6$}} 
\put(7,0) {\makebox(0,-1)[c]{\small $7$}} 
\put(8,0) {\makebox(0,-1)[c]{\small $8$}} 
\put(9,0) {\makebox(0,-1)[c]{\small $9$}} 
\put(10,0){\makebox(0,-1)[c]{\small $10$}}
\put(11,0){\makebox(0,-1)[c]{\small $11$}}
\put(12,0){\makebox(0,-1)[c]{\small $12$}}
\put(13,0){\makebox(0,-1)[c]{\small $13$}}
\put(14,0){\makebox(0,-1)[c]{\small $14$}}
\end{picture}}
\end{center}
\caption{A graphical representation of the permutation \usebox{\perm}
  corresponding to the PDSAW in Figure~\ref{fig:PDSAW}.b.}
  \label{fig:Permutation}
\end{figure}

Analogous to Theorem~\ref{thm:pdsaws} we can prove:

\begin{thm}\label{thm:pdsaws0}
  PDSAWs that stay below the $x$-axis, with $n$ east steps,
  \begin{itemize}
  \item $N$ north steps and
  \item length of the last descent equal to $M$
  \end{itemize}
  are in bijection with permutations of $[n]$ with
  \begin{itemize}
  \item $N$ nestings and
  \item $1$ being mapped to $M$.
  \end{itemize}
  Moreover, factors are preserved, the last factor of the PDSAW being mapped to
  the first factor of the permutation.
\end{thm}

\section{The bijection for PDSAWs in the symmetric wedge}\label{sec:symmetric}

In this section we exhibit a bijective proof of Theorem~\ref{thm:pdsaws}.  We
do so by taking a slight detour over certain weighted Dyck paths, known as \lq
histoires de Hermite\rq.

\subsection{A bijection between matchings and weighted Dyck paths}
\label{sec:KasraouiZeng}

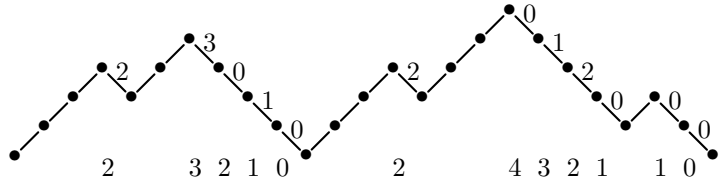
\begin{figure}
  \begin{equation*}
    \begin{xy}<11pt,0pt>:
      (4,-0.3)*{2},(7,-0.3)*{3},(8,-0.3)*{2},(9,-0.3)*{1},(10,-0.3)*{0},
      (14,-0.3)*{2},(18,-0.3)*{4},(19,-0.3)*{3},(20,-0.3)*{2},(21,-0.3)*{1},
      (23,-0.3)*{1},(24,-0.3)*{0},
      (4.5,3)*{2},(7.5,4)*{3},(8.5,3)*{0},(9.5,2)*{1},(10.5,1)*{0},
      (14.5,3)*{2},(18.5,5)*{0},(19.5,4)*{1},(20.5,3)*{2},(21.5,2)*{0},
      (23.5,2)*{0},(24.5,1)*{0},
      @i @={(1,1),(1,1),(1,1),(1,-1),(1,1),(1,1),(1,-1),(1,-1),(1,-1),(1,-1),(1,1),(1,1),(1,1),
        (1,-1),(1,1),(1,1),(1,1),(1,-1),(1,-1),(1,-1),
        (1,-1),(1,1),(1,-1),(1,-1)},
      (0.8,0.1)*{\bullet}="prev", @@{;p+"prev";"prev";**[|(2)]@{-}="prev"*{\bullet}}
    \end{xy}
  \end{equation*}
  \caption{A Dyck path with $12$ south-east steps corresponding to the PDSAW in
    Figure~\ref{fig:PDSAW}.a.  The heights of the south-east steps are written
    along the $x$-axis.}
\label{fig:Dyck}
\end{figure}

For convenience, we introduce a group of objects which are known to be in
bijection with matchings (see, for example, the article by Anisse Kasraoui and
Jiang Zeng~\cite{KasraouiZeng2006}), namely Dyck paths with weights on the
south-east steps.  A \Dfn{Dyck path} is a path starting at the origin, taking
north-east and south-east steps, returning to the $x$-axis but never going
below it.  The \Dfn{height} of a step is the $y$-coordinate of the point where
it ends, and we allow a non-negative weight on each south-east step, at most as
big as its height.  These objects are also known as \lq histoires de
Hermite\rq, because of their connection to the Hermite orthogonal polynomials.
An example of such a path is given in Figure~\ref{fig:Dyck}.  For brevity, we
will refer to these weighted Dyck paths always simply as \lq Dyck paths\rq.

To be able to keep track of all the statistics mentioned in
Theorem~\ref{thm:pdsaws}, we need to describe their meaning also on Dyck paths:
\begin{dfn}
  The \Dfn{total weight} of a Dyck path is the sum of the weights of its
  south-east steps.  The \Dfn{complementary weight} of a step is the difference
  of its height and its weight, and the \Dfn{complementary weight} of a Dyck
  path is the sum of the complementary weights of its steps.

  A \Dfn{factor} of a Dyck path is a subpath that starts and ends at the
  $x$-axis, but does not return to the $x$-axis otherwise.
\end{dfn}

The following are trivial consequences of the bijection exhibited in the
article by Anisse Kasraoui and Jiang Zeng~\cite{KasraouiZeng2006}, which is a
variant of bijections in Philippe Flajolet's and Xavier Viennot's
articles~\cite{MR592851,Viennot1984}:
\begin{prop}
  Matchings of the set $[2n]$ with
  \begin{itemize}
  \item $N$ nestings,
  \item $C$ crossings, and
  \item $1$ being matched with $M$,
  \end{itemize}
  are in bijection with Dyck paths with $n$ south-east steps,
  \begin{itemize}
  \item total weight $N$,
  \item complementary weight $C$ and
  \item $M$ being the position of the first south-east step with weight zero.
  \end{itemize}
  Moreover, factors are preserved.
\end{prop}

(The \Dfn{position} of the first step in a Dyck path is one, the second step
has position two, etc.)

\subsection{A bijection between weighted Dyck paths and PDSAWs}
\label{sec:bijection}

In this section we present a bijection between weighted Dyck paths and PDSAWs,
thus proving Theorem~\ref{thm:pdsaws}.  As a side remark, for PDSAWs without
north steps the transformation is particularly simple: it consists of rotating
the PDSAW counterclockwise by $45$\textdegree\ and reflecting the result at a
vertical line, giving all south-east steps weight zero.

Let $P$ be a PDSAW, given by the $y$-coordinates of its east steps.  We proceed
recursively: there is only one PDSAW with a single east step, so we assume that
$P$ has $n>1$ east steps.  Let $P^\prime$ be obtained from $P$ by removing its
last east step.

If the length of the last descent of $P$ is minimal, i.e., one, we map $P$ to
the Dyck path obtained by prepending
\begin{xy}<8pt,0pt>:
  (2,0.8)*{\scriptstyle 0},
  @i @={(1,1),(1,-1)},
  (0.1,0)*{\scriptstyle\bullet}="prev", 
  @@{;p+"prev";"prev";**[|(1.3)]@{-}="prev"*{\scriptstyle\bullet}}
\end{xy}
to the Dyck path $D^\prime$ corresponding to $P^\prime$.

Otherwise, let $\hat P$ be the PDSAW obtained from $P$ by lowering the last
east step by one and let $\hat D$ be the corresponding Dyck path.  Let $M(\hat
D)$ be the position of the first south-east step in $\hat D$ with weight zero.
We consider the step immediately before and the step immediately after, and
produce the Dyck path $D$ corresponding to $P$ according to the following five
local transformation rules:

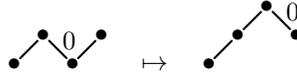
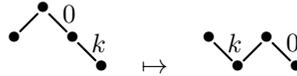
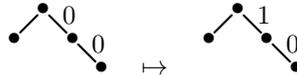
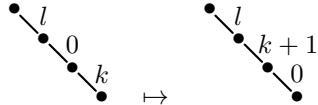
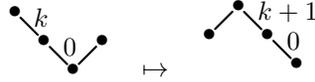
\renewcommand{\labelenumi}{(\Roman{enumi})}
\begin{enumerate}
\item


  \begin{equation*}
    \begin{xy}<11pt,0pt>:
      (2,0.8)*{0},
      @i @={(1,1),(1,-1),(1,1)},
      (0.1,0)*{\bullet}="prev", @@{;p+"prev";"prev";**[|(2)]@{-}="prev"*{\bullet}}
    \end{xy}\quad\mapsto\quad
    \begin{xy}<11pt,0pt>:
      (3,1.8)*{0},
      @i @={(1,1),(1,1),(1,-1)},
      (0.1,0)*{\bullet}="prev", @@{;p+"prev";"prev";**[|(2)]@{-}="prev"*{\bullet}}
    \end{xy}
  \end{equation*}

\item for $k>0$


  \begin{equation*}
    \begin{xy}<11pt,0pt>:
      (2,1.8)*{0},(3,0.8)*{k},
      @i @={(1,1),(1,-1),(1,-1)},
      (0.1,1)*{\bullet}="prev", @@{;p+"prev";"prev";**[|(2)]@{-}="prev"*{\bullet}}
    \end{xy}\quad\mapsto\quad
    \begin{xy}<11pt,0pt>:
      (1,0.8)*{k},(3,0.8)*{0},
      @i @={(1,-1),(1,1),(1,-1)},
      (0.1,1)*{\bullet}="prev", @@{;p+"prev";"prev";**[|(2)]@{-}="prev"*{\bullet}}
    \end{xy}
  \end{equation*}

\item 


  \begin{equation*}
    \begin{xy}<11pt,0pt>:
      (2,1.8)*{0},(3,0.8)*{0},
      @i @={(1,1),(1,-1),(1,-1)},
      (0.1,1)*{\bullet}="prev", @@{;p+"prev";"prev";**[|(2)]@{-}="prev"*{\bullet}}
    \end{xy}\quad\mapsto\quad
    \begin{xy}<11pt,0pt>:
      (2,1.8)*{1},(3,0.8)*{0},
      @i @={(1,1),(1,-1),(1,-1)},
      (0.1,1)*{\bullet}="prev", @@{;p+"prev";"prev";**[|(2)]@{-}="prev"*{\bullet}}
    \end{xy}
  \end{equation*}

\item


  \begin{equation*}
    \begin{xy}<11pt,0pt>:
      (1,2.8)*{l},(2,1.8)*{0},(3,0.8)*{k},
      @i @={(1,-1),(1,-1),(1,-1)},
      (0,3)*{\bullet}="prev", @@{;p+"prev";"prev";**[|(2)]@{-}="prev"*{\bullet}}
    \end{xy}\quad\mapsto\quad
    \begin{xy}<11pt,0pt>:
      (1,2.8)*{l},(2.7,1.8)*{k+1},(3,0.8)*{0},
      @i @={(1,-1),(1,-1),(1,-1)},
      (0,3)*{\bullet}="prev", @@{;p+"prev";"prev";**[|(2)]@{-}="prev"*{\bullet}}
    \end{xy}
  \end{equation*}

\item 
%

  \begin{equation*}
    \begin{xy}<11pt,0pt>:
      (1,1.8)*{k},(2,0.8)*{0},
      @i @={(1,-1),(1,-1),(1,1)},
      (0.1,2)*{\bullet}="prev",
      @@{;p+"prev";"prev";**[|(2)]@{-}="prev"*{\bullet}}
    \end{xy}\quad\mapsto\quad
    \begin{xy}<11pt,0pt>:
      (2.8,2.8)*{k+1},(3,1.8)*{0},
      @i @={(1,1),(1,-1),(1,-1)},
      (0.1,2)*{\bullet}="prev", @@{;p+"prev";"prev";**[|(2)]@{-}="prev"*{\bullet}}
    \end{xy}
  \end{equation*}
\end{enumerate}

We observe that for a given Dyck path, exactly one of these rules applies, and
each of them is invertible.  Furthermore, the rules imply that $M(D)=M(\hat
D)+1$.  It remains to check that also the other statistics in
Theorem~\ref{thm:pdsaws} are correctly mapped.

\begin{figure}
\begin{center}
{\setlength{\unitlength}{4mm}
\begin{picture}(17,8.6)(-8.5,-4.1)
\put(-0.85,-0.25){(III)}
\put(5.7,1.8){(IV)} 
\put(6.5,3.7){\circle{1.8}}       
\put(6.65,4.6){\vector(1,0){0}}    
\put(5.9,-2.3){(V)}
\put(6.5,-3.7){\circle{1.8}}      
\put(6.65,-4.6){\vector(1,0){0}}   
\qbezier(5.5,2)(4,0)(5.5,-2) 
\put(4.75,0){\vector(0,1){0}}
\qbezier(7.5,2)(9,0)(7.5,-2) 
\put(8.25,0){\vector(0,-1){0}}
\put(-7,1.8){(I)}
\put(-6.5,3.7){\circle{1.8}}
\put(-6.35,4.6){\vector(1,0){0}}
\put(-7.1,-2.3){(II)}
\put(-6.5,-3.7){\circle{1.8}}      
\put(-6.35,-4.6){\vector(1,0){0}}   
\qbezier(-5.5,2)(-4,0)(-5.5,-2)
\put(-4.75,0){\vector(0,1){0}}
\qbezier(-7.5,2)(-9,0)(-7.5,-2) 
\put(-8.25,0){\vector(0,-1){0}}
\put(-5,2){\vector(3,-1){4.5}}
\put(-5,-2){\vector(3,1){4.5}}
\put(0.5,0.5){\vector(3,1){4.5}}
\put(0.5,-0.5){\vector(3,-1){4.5}}
\end{picture}}
\end{center}
\caption{The automaton corresponding to the local rules from
  Section~\ref{sec:bijection}.}
  \label{fig:Automaton}
\end{figure}
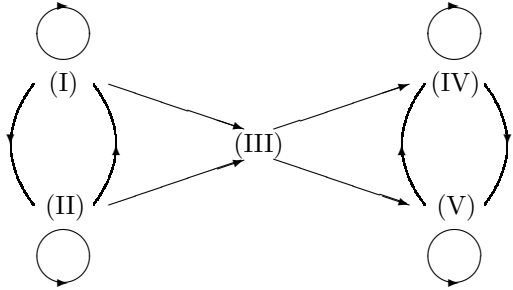

First, we remark that the five local transformation rules form an automaton as
in Figure~\ref{fig:Automaton}: transformation (I) can only be followed by one
of the transformations (I), (II) or (III), and so on.

We now show that the last two steps of the PDSAW have the same $y$-coordinate
precisely when transformation (III) applies, i.e., the first south-east step
with weight zero is preceded by a north-east step and followed by another
south-east step with weight zero:

Suppose that the last descent of the PDSAW $P^\prime$ has length $M-1$.  By
induction, $M$ is the position of the first south-east step of $D^\prime$ with
weight zero.  To obtain $D$, we first prepend
\begin{xy}<8pt,0pt>:
  (2,0.8)*{\scriptstyle 0},
  @i @={(1,1),(1,-1)},
  (0.1,0)*{\scriptstyle\bullet}="prev", 
  @@{;p+"prev";"prev";**[|(2)]@{-}="prev"*{\scriptstyle\bullet}}
\end{xy}
to $D^\prime$.  In this new path, the second south-east step with weight zero
is now at position $M+2$.  Then, if the last two east steps of $P$ have the
same $y$-coordinate, we must apply exactly $M-1$ times some local
transformation.  This will move the first south-east step having weight zero
from position $2$ to position $M+1$, i.e., just before the second south-east
step with weight zero.  We conclude that we only applied transformations of
type (I) and (II), and end up in a configuration as described in the preceding
paragraph.

Since (I) and (II) preserve the total weight of the Dyck path, and (III), (IV)
and (V) increase it by one, the bijection indeed transforms the number of north
steps into the total weight of the Dyck path.

A simple computation reveals that transformations (I) and (V) increase the
complementary weight by one, while the other transformations decrease it by
one.  Increasing the $y$-coordinate of the final east step of a PDSAW by one
also increases the area enclosed by the path and the line $y=-x$ by one.  Thus,
we find that the parity of the area and of the complementary weight coincides.

It remains to show that the bijection preserves prime factors, i.e., factors
that do not contain a smaller factor.  To start with, we observe that
transformations (II), (III) and (IV) preserve prime factors, while (I) and (V)
merge the first two prime factors, given that the middle step in the preimage
has height zero.  

The preimage of transformation (I), with middle step having height zero
corresponds exactly to the situation where the last east step of the PDSAW is
minimal.  The preimage must occur at the very beginning of the Dyck path, since
a step preceding it would have to be a south-east step at height zero, and
therefore necessarily of weight zero.

Finally, the preimage of transformation (V), with middle step having height
zero corresponds exactly to the situation where the last east step of the last
prime factor of the PDSAW is maximal, i.e., making it higher would merge the
last two prime factors.  This is because no south-east step before the middle
step in the preimage can have weight zero.  Therefore, the middle step of the
preimage marks the end of the first prime factor of the Dyck path, and the
length of this prime factor minus one is also the length of the last descent in
the PDSAW.

\section{The bijection for the asymmetric wedge}

The aim of this section is to present a construction for PDSAWs in the
asymmetric wedge analogous to the one presented in the preceding section.  The
r\^ole played by weighted Dyck paths in Section~\ref{sec:symmetric} is now
taken by weighted Motzkin paths, or \lq histoires de Laguerre\rq.  Note that
the bijection presented in the last section does \emph{not} seem to restrict
nicely to PDSAWs in the asymmetric wedge.  At least, we were unable to find a
good characterisation of those weighted Dyck paths that correspond to these
PDSAWs.

\subsection{A bijection between permutations and weighted bicoloured Motzkin
  Paths}
\label{sec:permutations}

We employ a slight variation of a bijection due to Dominique Foata and Doron
Zeilberger that maps permutations to weighted bicoloured Motzkin paths.

\begin{figure}
  \begin{equation*}
    \begin{xy}<11pt,0pt>:
      (2.3,1.7)*{1},(4.3,2.7)*{1},(5.5,2)*{1},(12.1,2)*{1},
      @i @={(1,1),(1,0),(1,1)},
      (0.8,0.1)*{\bullet}="prev",
      @@{;p+"prev";"prev";**[|(2)]@{-}="prev"*{\bullet}},
      (1,0);p+"prev";"prev";**[|(2)]@{=}="prev"*{\bullet},
      (1,-1);p+"prev";"prev";**[|(2)]@{-}="prev"*{\bullet},
      (1,0);p+"prev";"prev";**[|(2)]@{=}="prev"*{\bullet},
      @i @={(1,-1),(1,0),(1,1)},
      @@{;p+"prev";"prev";**[|(2)]@{-}="prev"*{\bullet}},
      (1,0);p+"prev";"prev";**[|(2)]@{=}="prev"*{\bullet},
      (1,0);p+"prev";"prev";**[|(2)]@{=}="prev"*{\bullet},
      @i @={(1,1),(1,-1),(1,-1)},
      @@{;p+"prev";"prev";**[|(2)]@{-}="prev"*{\bullet}},
    \end{xy}
  \end{equation*}
  \caption{A Motzkin path with $14$ steps corresponding to the PDSAW in
    Figure~\ref{fig:PDSAW}.b.  Only non-zero weights are indicated.}
\label{fig:Motzkin}
\end{figure}
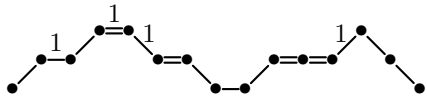

A \Dfn{bicoloured Motzkin path} is a path starting at the origin, taking
north-east, south-east, east and coloured east steps, returning to the $x$-axis
but never going below it.  The \Dfn{height} of a step is the $y$-coordinate of
the point where it ends, and we allow a non-negative weight on every step as
follows:
\begin{itemize}
\item[] south-east steps and east steps have weight at most as big as their
  height, and
\item[] north-east steps and coloured east steps have weight less than their
  height.
\end{itemize}
Note that this implies that coloured east steps must have height greater than
zero.  Such paths are also referred to as \lq histoires de Laguerre\rq.  An
example of such a path is given in Figure~\ref{fig:Motzkin}.  For brevity, we
will refer to weighted Motzkin paths always simply as \lq Motzkin paths\rq.

The \Dfn{total weight} of a Motzkin path is the sum of the weights of its
steps.  A \Dfn{factor} of a Motzkin path is a subpath that starts and ends at
the $x$-axis.  Thus, an east step of height zero constitutes a prime factor.

Following Sylvie Corteel~\cite{Corteel2007}, an arc $\big(i,\sigma(i)\big)$ of
a permutation is \Dfn{nested} by an arc $\big(j,\sigma(j)\big)$, if
$j<i\leq\sigma(i)<\sigma(j)$ or $j>i>\sigma(i)>\sigma(j)$.  Two arcs
$\big(i,\sigma(i)\big)$ and $\big(j,\sigma(j)\big)$ \emph{cross} if
$i<j\leq\sigma(i)<\sigma(j)$ (the middle inequality is weak!) or
$\sigma(j)<\sigma(i)<j<i$.

\begin{prop}
  Permutations of $[n]$ with
  \begin{itemize}
  \item $N$ nestings,
  \item $C$ crossings, and
  \item $1$ being mapped to $M$
  \end{itemize}
  are in bijection with Motzkin paths with $n$ steps,
  \begin{itemize}
  \item total weight $N$, 
  \item complementary weight $C$, and
  \item $M$ being the position of the first east or south-east step with weight
    $0$.
  \end{itemize}
  Moreover, factors are preserved.
\end{prop}

The path is created as follows: let $\sigma$ be a permutation of $[n]$, then
the $i$\textsuperscript{th} step is
\begin{itemize}
\item[] north-east, if $i<\min\left(\sigma(i),\sigma^{-1}(i)\right)$,
\item[] south-east, if $i>\max\left(\sigma(i),\sigma^{-1}(i)\right)$,
\item[] east, if $\sigma^{-1}(i)\leq i\leq\sigma(i)$, and
\item[] coloured east, if $\sigma(i)<i<\sigma^{-1}(i)$.
\end{itemize}
The weight of the $i$\textsuperscript{th} step is the number of arcs nesting
$\left(\sigma^{-1}(i), i\right)$.  We remark that this definition differs
slightly from Sylvie Corteel's in \cite{Corteel2007}, where the weight was
taken to be the number of arcs nesting $\left(i,\sigma(i)\right)$ instead.

Let us show that $\sigma(1)=M$ entails that the $M$\textsuperscript{th} step is
the first east or south-east step with weight zero: since it is impossible to
nest the arc $\big(1,\sigma(1)\big)$, the $M$\textsuperscript{th} step
certainly has weight zero.  Since either $M>\max\left(\sigma(M),
  \sigma^{-1}(M)=1\right)$ or $\sigma^{-1}(M)=1\leq M\leq\sigma(M)$, it must be
south-east or east.  Finally, suppose that the $i$\textsuperscript{th} step,
with $i<M$, is east or south-east. We then must have $\sigma^{-1}(i)\leq i$, so
it is nested by $\big(1,\sigma(1)\big)$ and has therefore non-zero weight.

To prove that prime factors are preserved, we remark that a left factor of a
permutation is a permutation by itself and therefore mapped to a Motzkin path,
and vice verse.

\subsection{A bijection between weighted Motzkin paths and PDSAWs below the
  $x$-axis}
\label{sec:bijection0}

For brevity, we refer in this section to PDSAWs that remain below the $x$-axis
simply as PDSAWs.  The bijection we are about to describe is very similar to
the one in Section~\ref{sec:bijection}, so we allow ourselves to keep the
description shorter.

Let $P$ be a PDSAW with $n>1$ east steps and let $P^\prime$ be obtained from
$P$ by removing its last east step.  If the length of the last descent of $P$
is minimal, we map it to the Motzkin path obtained by prepending an east step
to the Motzkin path $M^\prime$ corresponding to $P^\prime$.

Otherwise, let $\hat P$ be the PDSAW obtained from $P$ by lowering the last
east step by one and let $\hat M$ be the corresponding Motzkin path.  We
consider the first east or south-east step in $\hat M$ that has weight zero,
along with the step immediately after it. Then we produce $M$ corresponding to
$P$ according to the following local transformation rules:

\renewcommand{\labelenumi}{(\Roman{enumi})}
\begin{enumerate}
\item if the first step of the pair is an east step, and the second is not an
  east or south-east step with weight zero, exchange the two steps;
\item if the first step of the pair is an east step, and the second is an east
  step with weight zero, replace them with a north-east step followed by a
  south-east step, both with weight zero;
\item if the first step of the pair is an east step, and the second is a
  south-east step with weight zero, replace them with a coloured east step
  followed by a south-east step, both with weight zero;
\item otherwise, that is, if the first step of the pair is a south-east step,
  increase the weight of the second step by one and exchange the two steps and
  their weights.
\end{enumerate}

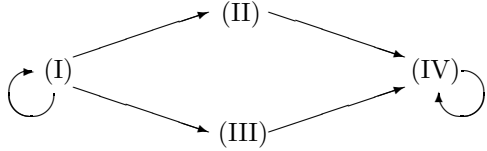
\begin{figure}
  \centering
{\setlength{\unitlength}{4mm}
\begin{picture}(12,4)(-6,-2)
\put(-6.5,-0.2){(I)}
\put(-6.9,-0.725){\oval(1.5,1.5)[b]\oval(1.5,1.5)[l]}
\put(-6.8,0){\vector(1,0){0}}
\put(-5.5,0.5){\vector(3,1){4.5}}
\put(-5.5,-0.5){\vector(3,-1){4.5}}
\put(1,2){\vector(3,-1){4.5}}
\put(1,-2){\vector(3,1){4.5}}
\put(-0.6,1.7){(II)}
\put(-0.7,-2.3){(III)}
\put(5.7,-0.2){(IV)}
\put(7.4,-0.725){\oval(1.5,1.5)[b]\oval(1.5,1.5)[r]}
\put(6.65,-0.65){\vector(0,1){0}}
\end{picture}}
\caption{The automaton corresponding to the local rules from
  Section~\ref{sec:bijection0}.}
  \label{fig:Automaton0}
\end{figure}

Again, these rules form an automaton, which is depicted in
Figure~\ref{fig:Automaton0}.  Note however that this automaton differs
structurally from the one in described in Section~\ref{sec:bijection}: here,
one of transformations (II) or (III) applies if the final east step of the
PDSAW is just below its next-to-last east step.  I.e., there are two transient
states, and the transition occurs \emph{before} the final two east steps of the
PDSAW have the same $y$-coordinate.  Since we enter state (IV) when the final
east step of $P$ is at the same height as the next-to last east step of $P$,
and each transformation in state (IV) increase the weight by one, every north
step between the last two east steps of $P$ contributes one to the total weight
of the resulting Motzkin path.

We have to remark that the bijection in this section does not transform the
number of crossings of the permutation into a meaningful statistic on PDSAWs.
This is a little disappointing, given the other parallels between
Theorems~\ref{thm:pdsaws} and \ref{thm:pdsaws0}.

\subsection{Philippe Nadeau's bijection}
\label{sec:nadeau}

In this section we present Philippe Nadeau's bijection, which was already
mentioned in the introduction. The theorem it proves is almost identical to
Theorem~\ref{thm:pdsaws0}:
\begin{thm}\label{thm:pdsaws0Philippe}
  PDSAWs that stay below the $x$-axis, with $n$ east steps,
  \begin{itemize}
  \item $N$ north steps and
  \item length of the last descent equal to $M$
  \end{itemize}
  are in bijection with permutations of $[n]$ with
  \begin{itemize}
  \item $N$ occurrences of the pattern $31\mbox{-}2$ and
  \item $1$ being mapped to $M$.
  \end{itemize}
  Moreover, prime factors are preserved, the last prime factor of the PDSAW
  being mapped to the first factor of the permutation.
\end{thm}

Let $P$ be a PDSAW with $n$ east steps and let $R=\{1,2,\dots,n\}$.  We
construct the corresponding permutation $\sigma$ as follows: for $i$ in
$1,2,\dots, n$,
\begin{itemize}
\item[] let $h$ be $1$ minus the $y$-coordinate of the
  $(n-i+1)$\textsuperscript{st} east step of $P$,
\item[] let $\sigma(i)$ be the $h$-largest element in $R$, and
\item[] delete this element from $R$.
\end{itemize}

It is immediate that $1$ is mapped to $M$, since for $i=1$ we have $h=1-(M-n)$.
To prove that north steps are translated to occurrences of the pattern
$31\mbox{-}2$, note that for every east step that is $k$ units above the
preceding step we introduce $k$ occurrences of the pattern $31\mbox{-}2$ in
$\sigma$.  (Of course, this does not mean, that the images of the steps under
$\sigma$ would differ by $k$!)

Apart from being extremely simple, this bijection has another beautiful
property: it is identical to the composition of the bijection presented in the
previous section with the bijection given by Sylvie Corteel in
\cite{Corteel2007}, although this does not appear obvious at all.

We have to remark, unfortunately, that we were unable to find a similar
bijection between matchings and PDSAWs in the symmetric wedge.  In particular,
we could not find a \lq nice\rq\ definition of a pattern in matchings, such
that the number of occurrences thereof would correspond to the number of north
steps in PDSAWs in the symmetric wedge.

\section{The end of the story?}

A brief look at the literature will convince us that there remains a fair bit
of work to be done.  Most pressing is the question whether we can derive the
generating functions obtained by Janse van Rensburg, Thomas Prellberg and
Andrew Rechnitzer~\cite{RensburgPrellbergRechnitzer2007} for PDSAWs ending
\emph{anywhere}, according to the \emph{total} number of steps.  For example,
for the symmetric wedge they obtain, with
\begin{align*}
  P &= \sqrt{(1-t^2)(1-5t^2)}\\
  Q &= (1-3t^2-P)/2t.
  \intertext{the generating function}
  \sum_{n\geq0}c_n t^n &= \frac{1}{t(1-2t-t^2)}\left((1+t)t-(1-t^2-P)
    \sum_{n\geq0} (-1)^n t^{n^2} Q^n\right)\\
  &= 1+t+3t^2+5t^3+13t^4+\dots
\end{align*}
Although the bijections described here treat the length of the last descent of
the PDSAW nicely, we were unable to derive a formula for PDSAWs for a fixed
value of this statistic, or a generating function keeping track of it.

Apart from that, it might be worth exploring whether there is a generalisation
of PDSAWs that correspond to partitions.  Since the bijection of Anisse
Kasraoui and Jiang Zeng also applies to these more general objects, it would
not be too surprising to find such a generalisation.  Maybe this would involve
PDSAWs allowing diagonal steps.

Related to this, we would like to point out the link to continued fractions and
orthogonal polynomials.  Indeed, the machinery developed by Philippe
Flajolet~\cite{MR592851} and Xavier Viennot~\cite{Viennot1984} teaches us to
interpret the expression in Equation~\eqref{eq:Touchard} as the
$2n$\textsuperscript{th} moment of a $q$-analogue of the Hermite polynomials,
and the expression in Equation~\eqref{eq:Williams} as the
$n$\textsuperscript{th} moment of a $q$-analogue of the Laguerre polynomials.
The generating function $\sum_{n\geq0} M_{2n}x^n$ has the continued fraction
expansion
\begin{equation*}
  \cfrac{1}{1-\cfrac{[1]_q x}{1-\cfrac{[2]_q x}{1-\cfrac{[3]_q x}{\ddots}}}}
\end{equation*}
while the generating function $\sum_{n\geq0} P_n x^n$ has the expansion
\begin{equation*}
  \cfrac{1}{1-[1]_q x -\cfrac{[1]_q^2 x^2}{1-([1]_q+[2]_q)x-
      \cfrac{[2]_q^2 x^2}{1-([2]_q+[3]_q)x-\cfrac{[3]_q^2 x^2}{\ddots}}}}
\end{equation*}

Finally, we would like to explain why we chose to present the bijections in
terms of automatons.  The reason, although realised only after having found the
bijections, is that they are unique in the following sense: given that the
bijections should preserve the parameters as described in
Theorems~\ref{thm:pdsaws}, and there should be \emph{local} rules displacing
the first occurrence of a south-east step with weight zero to the right, the
rules in Section~\ref{sec:bijection} are already determined.  Of course, this
does not exclude the possibility of other bijective proofs, as demonstrated
very recently by Svetlana Poznanovi\'c~\cite{Poznanovic2008}.


\providecommand{\cocoa} {\mbox{\rm C\kern-.13em o\kern-.07em C\kern-.13em
  o\kern-.15em A}}
\providecommand{\bysame}{\leavevmode\hbox to3em{\hrulefill}\thinspace}
\providecommand{\href}[2]{#2}

\end{document}